\documentclass{article}

\pdfoutput=1 

\usepackage{amsmath,amsfonts,amssymb} 
\usepackage[utf8]{inputenc}
\usepackage[pdftex]{graphicx} 
\usepackage{multirow} 

\newcommand{\Expc}{\mathop{\textrm{E}}}

\newcommand{\bbar}{\rule[-2.5pt]{1.5pt}{10pt}}

\title{Runs of Ones in Binary Strings}
\author{Félix Balado and Guénolé Silvestre}

\date{School of Computer Science\\University College Dublin, Ireland}

\usepackage{hyperref}
\hypersetup{
  colorlinks=true,
  urlcolor=blue,
  breaklinks=true
}
\begin{document}

\maketitle

\section{Introduction}

A run of ones in a binary string is an uninterrupted sequence of ones
flanked on each side either by a zero or by the start/end of the
string (Mood's criterion~\cite{mood40:_distr_theor_runs}). In the
following, a run of ones will be simply referred to as a ``run''.

What is the total number of runs of length $i$ over all binary
$n$-strings, or equivalently, if we draw binary $n$-strings uniformly
at random, what is the expected number of runs of length $i$ in a
binary $n$-string?

This problem was previously solved by
\href{https://www.researchgate.net/publication/254031393_Energy-efficient_communication_Understanding_the_distribution_of_runs_in_binary_strings}{Sinha
  and Sinha}~\cite{sinha12:_energy} ---in fact, these authors also
solved a harder problem in~\cite{sinha09} from which the solution to
the problem addressed here can be produced.

Here we wish to show that the solution can be found in shorter and
simpler ways than in~\cite{sinha09} or~\cite{sinha12:_energy} (see
also~\cite{makri12:_count_runs}).  We give three different solutions:
two of them use elementary counting arguments (recursive and
combinatorial solutions), while the third one is probabilistic.

\section{Counting Runs Recursively}\label{sec:counting-runs}
Let $r_n(i)$ be the total number of runs of length $i$ over all binary
$n$-strings.

First of all, let us get some visual intuition. In the diagrams below,
for $n=2,3$ and $4$ we list all binary $n$-strings and, right
underneath, their runs ``spectra'' (i.e. the number of runs of lengths
$1$ to $n$ found in each particular $n$-string). On the right we show
$r_n(i)$.

\begin{itemize}
\item $n=2$:

\begin{center}
  \begin{tabular}{c cccc c}
    &0&0&1&1&\\
    &0&1&0&1&\\
    \hline
    $i$&&&&& $r_2(i)$\\
    \hline
    \multicolumn{1}{c|}{1}& 0&1&1&0&\multicolumn{1}{|c}{2}\\
    \multicolumn{1}{c|}{2}& 0&0&0&1&\multicolumn{1}{|c}{1}\\
  \end{tabular}
\end{center}

\item $n=3$:

\begin{center}
  \begin{tabular}{c cccccccc c}
    &0&0&0&0&1&1&1&1&\\
    &0&0&1&1&0&0&1&1&\\
    &0&1&0&1&0&1&0&1&\\
    \hline
    $i$&&&&&&&&& $r_3(i)$\\
    \hline
    \multicolumn{1}{c|}{1}&0&1&1&0&1&2&0&0&\multicolumn{1}{|c}{5}\\
    \multicolumn{1}{c|}{2}&0&0&0&1&0&0&1&0&\multicolumn{1}{|c}{2}\\
    \multicolumn{1}{c|}{3}&0&0&0&0&0&0&0&1&\multicolumn{1}{|c}{1}\\
\end{tabular}
\end{center}

\item $n=4$:
  
\begin{center}
  \begin{tabular}{c cccccccccccccccc c}
    &0&0&0&0&0&0&0&0&1&1&1&1&1&1&1&1&\\
    &0&0&0&0&1&1&1&1&0&0&0&0&1&1&1&1&\\
    &0&0&1&1&0&0&1&1&0&0&1&1&0&0&1&1&\\
    &0&1&0&1&0&1&0&1&0&1&0&1&0&1&0&1&\\
    \hline  
    $i$&&&&&&&&&&&&&&&&&$r_4(i)$\\
    \hline
    \multicolumn{1}{c|}{1}&0&1&1&0&1&2&0&0&1&2&2&1&0&1&0&0&\multicolumn{1}{|c}{12}\\
    \multicolumn{1}{c|}{2}&0&0&0&1&0&0&1&0&0&0&0&1&1&1&0&0&\multicolumn{1}{|c}{5}\\
    \multicolumn{1}{c|}{3}&0&0&0&0&0&0&0&1&0&0&0&0&0&0&1&0&\multicolumn{1}{|c}{2}\\
    \multicolumn{1}{c|}{4}&0&0&0&0&0&0&0&0&0&0&0&0&0&0&0&1&\multicolumn{1}{|c}{1}\\
  \end{tabular}
\end{center}

\end{itemize}

Assuming $n>1$, it is trivial to see that
\begin{align}
  r_n(n)&=1,\nonumber\\
  r_n(n-1)&=2.\label{eq:rnnm1}
\end{align}
This solves the problem for $n=2$. Assuming $n>2$, let us next see
how, for $1\le i<n-1$, $r_n(i)$ can be recursively obtained from
$r_{n-1}(i)$:

\begin{itemize}
\item On the one hand, consider the $n$-strings that start
  with a zero: these trivially contribute $r_{n-1}(i)$ runs to $r_n(i)$.
\item On the other hand, consider the $n$-strings that start with a
  one. For $1\le i \le n-1$, the $2^{n-i-1}$ $n$-strings that start
  with $i$~ones followed by at least one zero add $2^{n-i-1}$ runs to
  $r_{n-1}(i)$ if $i<n-1$, and subtract $2^{n-i-1}$ runs from
  $r_{n-1}(i-1)$ if $i>1$.

  Thus, the $n$-strings that start with a one contribute
  $r_{n-1}(i)+2^{n-i-1}-2^{n-i-2}=r_{n-1}(i)+2^{n-i-2}$ runs to $r_n(i)$.
  
\end{itemize}
So, collecting these two contributions we have that
\begin{equation}\label{eq:rec}
  r_{n}(i)= 2\,r_{n-1}(i)+2^{n-i-2}.
\end{equation}
Now, by using the expression above recursively $k$ times, with
$k< n-i$, we get
\begin{equation}\label{eq:rn_k}
  r_n(i)=2^k r_{n-k}(i)+ k\, 2^{n-i-2}.
\end{equation}
When $k=n-i-1$, we have that $r_{n-k}(i)=r_{i+1}(i)=2$ because
of~\eqref{eq:rnnm1}. Thus, inputting this value of $k$
in~\eqref{eq:rn_k} we get the following explicit expression for the
number of runs of ones of length $i$ over all binary $n$-strings:
\begin{equation}\label{eq:rn}
  r_n(i)=(n-i+3)\,2^{n-i-2},
\end{equation}
where $1\le i<n-1$. Incidentally,~\eqref{eq:rn} is also valid when
$i=n-1$, i.e. it includes~\eqref{eq:rnnm1}. Furthermore, we can see
from~\eqref{eq:rn} that $r_{n-1}(i-1)=r_n(i)$ for $n\ge 2$, and so
recurrence~\eqref{eq:rec} can alternatively be expressed as a
recurrence on~$i$ for a given $n$ (rather than as a recurrence on $n$ for a
given $i$) as
\begin{equation}
  r_{n}(i)= 2\,r_{n}(i+1)+2^{n-i-2}.\label{eq:rec_alt}
\end{equation}

Also, the total number of runs over all $n$-strings is
\begin{equation}
  \label{eq:total}
  t(n)=\sum_{i=1}^n r_n(i)=(n+1)\,2^{n-2}.
\end{equation} %

Recursion~\eqref{eq:rec} was previously given
in~\cite{sinha12:_energy}
and~\cite{makri12:_count_runs}. Also~\eqref{eq:rn}
and~\eqref{eq:total} were given in~\cite{sinha12:_energy}.

\section{Counting Runs Using Combinatorics}
\label{sec:counting-runs-using}
Let us now give a way to obtain~\eqref{eq:total} directly by using
combinatorics, i.e. without relying on~\eqref{eq:rn}.  Consider the
partitioning of an arbitrary $n$-string into $p$ nonempty substrings,
where $1\le p\le n$. As shown in the example below, in which the
$n$-string is represented by $n$ asterisks and the $p$ partitions by $p+1$
vertical bars, we can put that partition into a bijection with $p$
runs of ones from binary $n$-strings:
\begin{center}
  \begin{tabular}{@{}c@{}c@{}c@{}c@{}c@{}c@{}c@{}c@{}c@{}c@{}c@{}c@{}c@{}c@{}c@{}c@{}c@{}c@{}c@{}c@{}c@{}c@{}c@{}c@{}c@{}c@{}c@{}}
    $\mid$&$*$&$*$&$*$&$\mid$&$*$&$*$&$*$&$*$&$\mid$& $*$&$*$ &$\mid$& $*$&$*$&$*$&$*$&$*$ &$\mid$&$\cdots$ &$\mid$&$*$&$\mid$&$*$&$*$&$\mid$\\
    \hline
    &1&1&1&  &0&0&0&0 & & 1&1 & & 0&0&0&0&0 &&$\cdots$ && 0&& 1&1&\\
    &0&0&0&  &1&1&1&1 & & 0&0 & & 1&1&1&1&1 && $\cdots$ && 1&& 0&0&\\
  \end{tabular}
\end{center}

The number of ways in which we can partition an $n$-string into $p$
nonempty substrings is ${n-1\choose p-1}$.  As we have $p$ runs of
ones associated to each of the ${n-1\choose p-1}$ possibilities, the
total number of runs over all binary $n$-strings
is%
\begin{equation}
  \label{eq:total_alt}
  t(n)=\sum_{p=1}^{n} p \,{n-1\choose p-1}=(n+1)\,2^{n-2}.
\end{equation}

Next, let us exploit the parallelism between expressions~\eqref{eq:rn}
and~\eqref{eq:total} to write $r_n(i)$ in a way that
echoes~\eqref{eq:total_alt}. In order to do so, let us put
\eqref{eq:rn} as the addition of the following two terms:
\begin{equation}
  \label{eq:rn_expanded}
  r_n(i)=(n-i+1)\,2^{n-i-2}+2^{n-i-1}.
\end{equation}
We can now see that the first term can expanded using the same
summation as in~\eqref{eq:total_alt}. Interestingly, by using
$2^n=\sum_{k=0}^n {n\choose k}$, the second term can also be expanded
into a summation with the same limits and over the same binomial
coefficients as the first one, but without the $p$ multiplication
factors. Combining these two observations we can
write~\eqref{eq:rn_expanded} as follows:
\begin{align}
  r_n(i)&=\sum_{p=1}^{n-i}p\, {n-i-1\choose p-1}+\sum_{p=1}^{n-i}{n-i-1\choose p-1}\nonumber\\
        &=\sum_{p=1}^{n-i}(p+1)\, {n-i-1\choose p-1}.  \label{eq:rn_expanded_2}
\end{align}
Expression~\eqref{eq:rn_expanded_2} literally tells us that there is a
one-to-one correspondence between a partitioning of an arbitrary
$(n-i)$-string intro $p$ nonempty substrings and $p+1$ runs of length
$i$ from binary $n$-strings. This provides a simple combinatorial
approach to obtaining~\eqref{eq:rn}.

In order to see how the bijection implied by~\eqref{eq:rn_expanded_2}
works we just need to place a run of ones of length $i$ at each of the
$p+1$ possible positions in the partitioning (i.e. the $p-1$ divisions
plus the start and the end of the string) and then complete a binary
$n$-string with alternating runs of ones and zeros having lengths
determined by the partitioning.

For example, let us take $n=4$ and $i=1$ and explicitly illustrate the
correspondence between each of the possible partitionings of an
arbitrary $3$-string ``$***$'' into~$p$~nonempty substrings and $p+1$
runs of length $1$ in binary $4$-strings:
\begin{itemize}
\item $p=1$ $\to$ ${2\choose 0}=1$ partitioning associated with $2$
  runs of length $1$

  \begin{flushleft}
    \begin{tabular}{c@{}c@{}c@{}c@{}c l}
      $\mid$&$*$&$*$&$*$&$\mid$ & \\
      \hline
      $\bbar$&$*$&$*$&$*$&$\mid$ & $\mathbf{1}000$\\
      $\mid$&$*$&$*$&$*$&$\bbar$ & $000\mathbf{1}$ 
    \end{tabular}
  \end{flushleft}
  
\item $p=2$ $\to$ ${2\choose 1}=2$ partitionings, each associated with
  $3$ runs of length $1$:

  \begin{flushleft}
    \begin{tabular}{c@{}c@{}c@{}c@{}c@{}c l}
      $\mid$&$*$&$\mid$&$*$&$*$&$\mid$ & \\
      \hline
      $\bbar$&$*$&$\mid$&$*$&$*$&$\mid$ & $\mathbf{1}011$\\
      $\mid$&$*$&$\bbar$&$*$&$*$&$\mid$ & $0\mathbf{1}00$\\
      $\mid$&$*$&$\mid$&$*$&$*$&$\bbar$ & $100\mathbf{1}$
    \end{tabular}
  \end{flushleft}

  \begin{flushleft}
    \begin{tabular}{c@{}c@{}c@{}c@{}c@{}c l}
      $\mid$&$*$&$*$&$\mid$&$*$&$\mid$ & \\
      \hline
      $\bbar$&$*$&$*$&$\mid$&$*$&$\mid$ & $\mathbf{1}001$\\
      $\mid$&$*$&$*$&$\bbar$&$*$&$\mid$ & $00\mathbf{1}0$\\
      $\mid$&$*$&$*$&$\mid$&$*$&$\bbar$ & $110\mathbf{1}$\\
    \end{tabular}
  \end{flushleft}

\item $p=3$ $\to$ ${2\choose 2}=1$ partitioning, associated with $4$
  runs of length $1$:

  \begin{flushleft}
    \begin{tabular}{c@{}c@{}c@{}c@{}c@{}c@{}c l}
      $\mid$&$*$&$\mid$&$*$&$\mid$&$*$&$\mid$ & \\
      \hline
      $\bbar$&$*$&$\mid$&$*$&$\mid$&$*$&$\mid$ & $\mathbf{1}010$\\
      $\mid$&$*$&$\bbar$&$*$&$\mid$&$*$&$\mid$ & $0\mathbf{1}01$\\
      $\mid$&$*$&$\mid$&$*$&$\bbar$&$*$&$\mid$ & $10\mathbf{1}0$\\
      $\mid$&$*$&$\mid$&$*$&$\mid$&$*$&$\bbar$ & $010\mathbf{1}$\\
    \end{tabular}
  \end{flushleft}
\end{itemize}

If $n=4$ and $i=2$ we look at the partitionings of ``$* *$'' into $p$
nonempty substrings:

\begin{itemize}
\item $p=1$ $\to$ ${1\choose 0}=1$ partitioning associated with $2$
  runs of length $2$
  \begin{flushleft}
    \begin{tabular}{c@{}c@{}c@{}c l}
      $\mid$&$*$&$*$&$\mid$ & \\
      \hline
      $\bbar$&$*$&$*$&$\mid$ & $\mathbf{11}00$\\
      $\mid$&$*$&$*$&$\bbar$ & $00\mathbf{11}$ 
    \end{tabular}
  \end{flushleft}
\item $p=2$ $\to$ ${1\choose 1}=1$ partitioning associated with $3$
  runs of length $2$

  \begin{flushleft}
    \begin{tabular}{c@{}c@{}c@{}c@{}c l}
      $\mid$&$*$&$\mid$&$*$&$\mid$ & \\
      \hline
      $\bbar$&$*$&$\mid$&$*$&$\mid$ & $\mathbf{11}01$\\
      $\mid$&$*$&$\bbar$&$*$&$\mid$ & $0\mathbf{11}0$\\
      $\mid$&$*$&$\mid$&$*$&$\bbar$ & $10\mathbf{11}$\\
    \end{tabular}
  \end{flushleft}
\end{itemize}

Finally, if we let the sum in~\eqref{eq:rn_expanded_2} start at
$p=0$ rather than $p=1$, then~\eqref{eq:rn_expanded_2} also includes
the case\footnote{Using ${-1\choose -1}=1$~\cite{kronenburg11} and
  ${t\choose -1}=0$ for all nonnegative integers~$t$.} $i=n$,
whereas~\eqref{eq:rn} is only valid for $i\le n-1$.

\section{Counting Runs Probabilistically}
\label{sec:count-runs-prob}
Assume that the binary $n$-strings are generated uniformly at random,
and define~$n$ indicator random variables $S_1,\ldots,S_n$ such that
$S_k=1$ if a run of ones starts at the $k$-th position of the string
and $0$ otherwise. There are two cases for the probability of $S_k=1$:
\begin{equation}
  \label{eq:pr_sk}
  \Pr(S_k=1)=\left\{
    \begin{array}{l l}
      \frac{1}{2}&\text{ if } k=1\\
      \frac{1}{4}&\text{ if } 2\le k\le n
    \end{array}\right.,
\end{equation}
as if a run starts at position $k>1$ then the previous bit must be
zero. The total number of runs is thus given by r.v. $T=\sum_{k=1}^n S_k$,
whose expectation $\Expc(T)=\sum_{k=1}^n \Expc(S_k)=\sum_{k=1}^n \Pr(S_k=1)$ is
\begin{align}
  \label{eq:exp_t}
  \Expc(T)&=(n+1)\,2^{-2}.
\end{align}
The rationale above was given by Alex Proscurin
in~\href{https://math.stackexchange.com/questions/1048541/expected-value-of-a-run-of-a-random-bitstring}{Expected value of a run of a random bitstring}.  As we can
see,~\eqref{eq:exp_t} is the normalisation of~\eqref{eq:total}, i.e.
$\Expc(T)=t(n)/2^n$.

Let us extend the procedure above and define indicator r.v.'s
$R_k^{(i)}$ such that $\Pr(R_k^{(i)}=1)$ is the probability that a run
of length $i$ starts at position $k$.  We have that
\begin{equation}
\Pr(R_k^{(i)}=1)=\Pr(R_k^{(i)}=1|S_k=1)\Pr(S_k=1),\label{eq:prk}
\end{equation}
because $\Pr(R_k^{(i)}=1|S_k=0)=0$. As for the conditional
probabilities in~\eqref{eq:prk}, there are three cases:
\begin{equation}
  \label{eq:cond_pr}
  \Pr(R_k^{(i)}=1|S_k=1)=\left\{
  \begin{array}{l l}
    2^{-i}& \text{ if }1\le k\le n-i\\
    2^{-(i-1)}& \text{ if }k= n-i+1\\
    0 & \text{ if } n-i+1<k\le n
  \end{array}\right..
\end{equation}
The first two cases are just the application of the geometric
distribution with parameter $1/2$, noting that in the second case the
run does not need to be finished by a zero. Now, the number of runs of
length $i$ is given by r.v. $T^{(i)}=\sum_{k=1}^n R_k^{(i)}$, whose
expectation $\Expc(T^{(i)})=\sum_{k=1}^n \Expc(R_k^{(i)})=\sum_{k=1}^n \Pr(R_k^{(i)}=1)$ is
\begin{equation}
  \label{eq:exp_ti}
  \Expc(T^{(i)})=(n-i+3)\,2^{-i-2}.
\end{equation}
Again, we have that~\eqref{eq:exp_ti} is the normalisation of~\eqref{eq:rn}, i.e. $\Expc(T^{(i)})=r_n(i)/2^n$.

The total fraction of runs of length $i$ is
\begin{equation}
  \label{eq:fni}
  f_n(i)=\frac{r_n(i)}{t(n)}=\frac{n-i+3}{n+1}\,2^{-i}
\end{equation}
for $1\le i\le n-1$, whereas $f_n(n)=1/t(n)$. Interestingly,
$f_n(2)=1/4$ for any $n>2$.  For large $n$,~\eqref{eq:fni} should
approximate the fraction of runs of length $i$ in one single binary
$n$-string drawn uniformly at random, in which case
$f_n(i)\approx 2^{-i}$. %

\section{Connections with Compositions, and Relationship to OEIS
  Sequences}
\label{sec:relat-sequ-a045623-1}

Sequence \href{https://oeis.org/A045623}{A045623} from
OEIS~\cite{oeis} (number of $1$'s in all compositions of $j+1$) is
defined by $a(j)=(j+3)\,2^{j-2}$ for $j\ge 1$ and
$a(0)=1$. Observing~\eqref{eq:rn}, $r_n(i)=a(n-i)$. 

Also, sequence \href{https://oeis.org/A001792}{A001792} (number of
parts in all compositions of $j + 1$) is defined
by $b(j) = (j+2)\, 2^{j-1}$. From~\eqref{eq:total}, $t(n)=b(n-1)$.

Similarly, the equivalent of expression~\eqref{eq:total_alt}
and expression~\eqref{eq:rn_expanded_2} are given in
\href{https://oeis.org/A045623}{A045623} and
\href{https://oeis.org/A001792}{A001792}, respectively.

These connections are due to the bijection between runs in binary
$n$-strings and parts of the compositions of $n$, which is easy to
observe when the compositions are put as the different ways in which
the sum $1+1+\cdots+1=n$ can be partitioned into nonzero parts (see
diagram at start of Section~\ref{sec:counting-runs-using}). Therefore,
the runs in all binary $n$-strings are in a one-to-one correspondence
with the parts in all compositions of $n$, which is $t(n)$ in both
cases. Similarly, the runs of length $1$ in all binary $n$-strings are
in a one-to-one correspondence with the $1$'s in all compositions of
$n$, which is $r_n(1)$ in both cases.

After noticing the connection between runs in binary strings and
compositions, we found that the gist of the counting argument that we
used in Section~\ref{sec:counting-runs} was previously given by Mike
Earnest in
\href{https://math.stackexchange.com/questions/3531833/number-of-parts-equal-to-1-in-all-compositions-of-n}{Number
  of parts equal to $1$ in all compositions of $n$}. Also, the
argument to get~\eqref{eq:total_alt} must be well known in the
compositions problem.

\section*{Acknowledgements}
\label{sec:acknowledgements}
Thanks to Kevin Ryde for pointing out the connection with A001729.

\bibliographystyle{unsrt}

\end{document}